\documentclass{amsart}

\usepackage[left=1.5in,top=1.0in,right=1.5in]{geometry}
\usepackage{amsthm}
\usepackage{amsmath}
\usepackage{amssymb,subcaption}
\usepackage{tikz}
\usepackage{euler}
\usepackage{url,hyperref}
\usepackage{colonequals}
\usepackage{mathrsfs}

\newtheorem{theorem}{Theorem}[section]
\newtheorem{proposition}[theorem]{Proposition}

\newtheorem{lemma}[theorem]{Lemma}

\newtheorem{lettertheorem}{Theorem}

\theoremstyle{definition}
\newtheorem{definition}[theorem]{Definition}

\newtheorem{example}[theorem]{Example}
\newtheorem{remark}[theorem]{Remark}

\numberwithin{equation}{section}

\def\N{{\mathbb N}}

\def\R{{\mathbb R}}

\def\A{{\mathcal A}}
\def\D{{\mathcal D}}
\def\V{{\mathcal V}}

\def\hU{\overline{U}}
\def\hu{\overline{\textbf{u}}}

\def\F{{\mathcal F}}
\def\rR{{\mathcal R}}
\def\G{{\mathcal G}}
\def\O{{\mathbf 0}}

\def\bd{\textbf{d}}

\def\conv{\operatorname{Conv}}
\def\gale{\operatorname{Gale}}
\def\cone{\operatorname{Cone}}

\newcommand{\Mink}{\mathtt{Mink}}
\newcommand{\TMink}{\mathtt{TMink}}
\newcommand{\DMink}{\mathtt{DMink}}

\newcommand{\relint}{\mathrm{relint}}

\newcommand{\type}{\mathbb{TP}}

\newcommand{\typecone}{\mathbb{TC}}

\definecolor{orange}{rgb}{0.898, 0.621, 0.0}
\definecolor{skyblue}{rgb}{0.336, 0.703, 0.910}
\definecolor{bluishgreen}{rgb}{0, 0.617, 0.449}
\definecolor{yellow}{rgb}{0.937, 0.890, 0.258}
\definecolor{blue}{rgb}{0, 0.445, 0.695}
\definecolor{red}{rgb}{0.832, 0.367, 0}
\definecolor{purple}{rgb}{0.797, 0.473, 0.652}

\makeatletter
\newcommand{\barefootnote}[1]{{%
  \let\@makefnmark\relax \let\@thefnmark\relax
  \makeatletter
  \@footnotetext{#1}}}
\makeatother 

\title{
	Minkowski summands of Cubes 
}
\author[F. Castillo]{Federico Castillo}
\author[J. Doolittle]{Joseph Doolittle}
\author[B. Goeckner]{Bennet Goeckner}
\author[M. S. Ross]{Michael S. Ross}
\author[L. Ying]{Li Ying}

\begin{document}
	\begin{abstract} 
		In pioneering works of Meyer and of McMullen in the early 1970s, the set of Minkowski summands of a polytope was shown to be a polyhedral cone called the type cone. Explicit computations of type cones are in general intractable. Nevertheless, we show that the type cone of the product of simplices is simplicial. This remarkably simple result derives from insights about rainbow point configurations and the work of McMullen.
	\end{abstract}
	\maketitle
	
	\section{Introduction}
	
	\barefootnote{2010 Mathematics Subject Classification: 52B05, 52B35} 
	
	A fundamental operation on polytopes is Minkowski addition. In this paper, we consider the reverse of this operation. Our motivating question is ``Given a polytope $P$, what can we say about the set of its Minkowski summands?''
	
	It is convenient to modify this question and consider the set $\TMink(P)$ of weak Minkowski summands, polytopes that are summands of some positive dilate of $P$, up to translation equivalence. It turns out that $\TMink(P)$ can be realized as a pointed polyhedral cone, which we refer to as the type cone of $P$.
	When $P$ is a rational polytope, the type cone is the nef cone of the associated toric variety associated to $P$, see \cite[Chapter~6]{toric}.
	We briefly describe three different parametrizations of the type cone.

	The starting point is a theorem of Shephard \cite[Section 15]{grunb} characterizing the weak Minkowski summands in terms of their support functions. In \cite{meyer}, Meyer used this connection to give a parametrization of $\TMink(P)$ using one parameter for each facet, so we refer to this construction as the facet parametrization.
	This parametrization has been used recently to compute type cones; see, e.g., \cite{coxeter,  nested, chapoton, pilaud_arnau}. Notably, the type cone of the regular permutohedron is the cone of submodular functions \cite{edmonds}.

	Shephard's aforementioned theorem provides another characterization: $Q$ is a weak Minkowski summand of \(P\) if and only if we can obtain $Q$ by moving the vertices of $P$ while preserving edge directions, also allowing contraction of edges to points. It follows that we can parametrize weak Minkowski summands by the edge lengths. This parametrization is called the \emph{edge deformation space} in \cite{faces} and is equal to the set of nonnegative 1-Minkowski weights. The set of $r$-Minkowski weights, as defined in \cite[Section 5]{simple} and further explored on \cite{weights}, is crucial for the understanding of McMullen's polytope algebra.
	
	Finally, in \cite{mcmullen} McMullen used a different description, using the support function as in \cite{meyer}, but expressing the whole set as an intersection of cones, one for each cofacet. McMullen calls these sets type cones, since the interior of such a cone parametrizes polytopes of a strong combinatorial type\footnote{McMullen \cite[Section 2]{weights} uses ``strongly isomorphic'' to refer to polytopes with the same normal fan.}. Abusing notation, in this paper we use the term type cone to refer to the \emph{closure} of what \cite{mcmullen} calls type cone. 
	
	We present basic facts about type cones using the edge parametrization in Section~\ref{sec:polygons}, even though we apply McMullen's method to obtain our main result. Each parametrization has shown its worth in different contexts. Facet parametrization has been successfully used on the regular permutohedron,  edge parametrization on polygons, and Gale intersections on cubes. But any mixing of these methods and polytopes has not been fruitful.
	
	Our main result is a description of type cones for $d$-cubes\footnote{By $d$-cubes, we mean any polytope combinatorially isomorphic to $[0,1]^d$.} and more generally for polytopes that are combinatorially isomorphic to products of simplices. Cubes can have quite nontrivial geometry. For instance, Klee and Minty in \cite{klee-minty} famously constructed cubes for which Dantzing's simplex method takes exponentially many steps. Also, for $d \geq 3$, there exist $d$-cubes for which each pair of opposite facets is orthogonal \cite{eg-models}. Surprisingly, cubes have elementary type cones. 
	
	\begin{lettertheorem}\label{thm:market}
		The type cone of any $d$-cube is a $d$-simplicial cone. More generally, any polytope that is combinatorially isomorphic to a product of simplices has a simplicial type cone.
	\end{lettertheorem}
	
	In contrast, Example \ref{ex:polygons} shows that the type cones of polygons are as general as possible within the dimension and facet count constraints.
	This gives an indication that the complexity of computing the type cone of a polytope cannot be easily determined from the complexity of the polytope itself.
	
	Recently Adiprasito, Kalmanovich, and Nevo proved that the realization space of the cube is contractible in \cite{karim}. Realization spaces parametrize the set of combinatorially isomorphic polytopes, whereas the interiors of type cones parametrize the set of polytopes with identical normal fans.
	
	The structure of realization spaces as semialgebraic sets can be arbitrarily complicated, and often one can only study their topology; see, e.g., \cite{mnev} and \cite{realizations}. 
	In contrast to realization spaces, type cones are topologically trivial, leaving their face structure as the primary object of interest. With the main result of \cite{karim} we find an interesting distinction. The realization space of both polygons and cubes are contractible, but the possible type cones have dramatically different combinatorics.

	When $P$ has a simplicial type cone, \cite[Corollary 1.11]{pilaud_arnau}	shows an explicit isomorphism between the type cone of \(P\) and the positive orthant. This isomorphism appears in the construction of Arkani-Hamed, Bai, He, and Yan \cite[Section 3.2]{scattering} of the \emph{kinematic} associahedron in the context of scattering amplitudes. Motivated by this connection with theoretical physics, Padrol, Palu, Pilaud, and Plamondon \cite{pilaud_arnau} analyzed the type cones of several families of polytopes to determine when they are simplicial. For instance, the only graph associahedra with simplicial type cones are the paths, i.e., the classical Loday's associahedra \cite{loday}. We remark that this is not true for other associahedra, in dimension $2$ associahedra are simply pentagons and by Example \ref{ex:polygons} their type cones can have 3, 4, or 5 facets.
	In another recent paper, Albertin, Pilaud, and Ritter in \cite{apr} classified which permutrees have simplicial type cones.

	From the results in \cite{apr} and \cite{pilaud_arnau}, it seems that having a simplicial type cone is a rare property. Moreover, these results depend on particular realizations, whereas our Theorem~\ref{thm:main} shows that \emph{all} realizations of products of simplices have simplicial type cones. Are these the only simple polytopes with this property?

	\section*{Acknowledgements}
	
	This work was completed in part at the 2019 Graduate Research Workshop in Combinatorics, which was supported in part by NSF grant \#1923238, NSA grant \#H98230-18-1-0017,  a generous award from the Combinatorics Foundation, and Simons Foundation Collaboration Grants \#426971 (to M. Ferrara) and \#315347 (to J. Martin). We thank Margaret Bayer and Tyrrell McAllister for discussions and encouragements in the early stages of the project.
	
	We are particularly grateful to Jeremy Martin, Isabella Novik, Arnau Padrol, Vincent Pilaud, and Raman Sanyal for many insightful discussions and comments on earlier drafts. We thank Jean-Philippe Labb\'e for computational help and Jes\'us De Loera for discussions about oriented matroids.  We thank the referees for their careful reading and helpful remarks.
	
	\section{Background}\label{sec:background}
	
	Let $\R^d$ be the $d$-dimensional Euclidean space with the usual inner product $\langle \cdot,\cdot \rangle:\R^d\times \R^d\to\R$. An \textbf{$n$-point configuration} is a $d\times n$ matrix $M_\A$. We think of it via its multiset of columns $\A=\{a_1,\dots, a_n\}\subset \R^d$. Abusing notation, we identify a point $a_i$ with its label $i$.
	A \textbf{face} of a point configuration \(\A\) is a subset \(S \subset \A\) such that there exists \(c\in \mathbb{R}^d\) with the following property: For any $y\in S$ and $x\in \A$ we have \(\langle c,x \rangle \leq \langle c, y\rangle\) with \(\langle c,x \rangle = \langle c, y\rangle\) if and only if $x\in S$. 
	We further include the empty set as a face of \(\A\), and note that the empty set and \(\A\) itself are called improper faces.
	The \textbf{dimension} of a face is the dimension of its affine hull.
	The set of all $k$-dimensional faces of $\A$ is denoted $\F_k(\A)$.  A \textbf{vertex} is a face of dimension \(0\), an \textbf{edge} is a face of dimension \(1\), and a \textbf{facet} is a face of codimension \(1\). A \textbf{coface} is a set of points \(S\subset\A\) such that \(\A \setminus S\) is a face, and a \textbf{cofacet} is the coface of a facet. If $\dim(\A)=d$, the vector $f(\A) \colonequals (f_0(\A),\dots,f_d(\A))$, where $f_k(\A) \colonequals |\F_k(\A)|$, is called the \textbf{f-vector} of \(\A\). The set of faces of \(\A\) forms a partially ordered set under inclusion called the \textbf{face lattice} $\F(\A)$.  Two point configurations $\A_1$ and $\A_2$ are said to be \textbf{combinatorially isomorphic} if their face lattices are isomorphic. 
	
	The \textbf{polytope} $P=P_\A$ obtained from a point configuration \(\A\)  is the convex hull $P_\A \colonequals \conv(\A)$.  Given a polytope $P$, we may treat it as a point configuration $\V(P)$ whose points are the vertices of \(P\) in some order. A cone is  \textbf{simplicial} if it has a linearly independent system of generators. Equivalently, a simplicial cone is the cone over a simplex.
	
	The polytope $\Delta_d:=\conv\{e_0,\dots,e_{d}\}\subset \R^{d+1}$, where the $e_i$ are the standard basis vectors, is called the \textbf{standard simplex} of dimension $d$. Its face lattice is the boolean lattice $B_{d+1}$ since every subset of the vertices forms a face. Any polytope combinatorially isomorphic to $\Delta_d$ is called a $d$-simplex, or simply a simplex if we do not specify dimension.
	
	Given a polytope $P\subset\R^d$ with $\dim P< d$, we can restrict to its affine hull, where it is full dimensional. Also, after some translation, any $d$-dimensional polytope $P$ in $\R^d$ contains the origin in the interior. In the present paper there is no harm in assuming that $P$ is full dimensional and contains the origin in the interior, in which case we define the \textbf{polar} polytope $P^\circ:=\{c \in \R^d : \langle c, x \rangle \leq 1 \text{ for all }x \in P\}$. On the level of face lattices, the face lattice of \(P^\circ\) is isomorphic to the face lattice of \(P\) with the order reversed. 
	
	By the Weyl--Minkowski Theorem \cite[Theorem 1.1]{ziegler}, a polytope $P$ can be alternatively described as the solution set to a finite system of linear inequalities, i.e., a \(d\)-dimensional polytope $P=\{x\in \R^d: Ux\leq z\}$ where $U$ is a $m\times d$ matrix and $z\in \R^m$. 
	If deleting any row of $U$ changes $P$, we call the system \textbf{irredundant} or \textbf{facet-defining}, since in this case each set $\{x\in P: \langle u_i,x\rangle = z_i\}$ defines a facet of $P$.

	\begin{remark}
		Any $d$-polytope $P$ with the origin in the interior can be presented as \begin{equation}\label{eq:system}
		P=\{x\in\R^d: Ux\leq 1\}
		\end{equation}
		for some matrix $U$ with $d$ columns and we allow the system to be redundant.  From the system we can read the polar polytope as \(P^{\circ}=
		\conv\{u: u\in \text{Rows}(U)\}\).
		
	\end{remark}
	\begin{definition}\label{def:polar}
		
		For our purposes we need a slightly more general notion of polarity.
		Given a system of inequalities of the form given in Equation \eqref{eq:system}, we define the $\D(P)$ as the point configuration of the row vectors of the system. The convex hull of $\D(P)$ is $P^\circ$. A point $u_i$ is a vertex if and only if $\langle u_i,x\rangle =1$ defines a facet of $P$.\end{definition}
	
	Let $Q\subset \R^c, R \subset \R^d$ be two polytopes. Their \textbf{(Cartesian) product} is
	\[
	Q\times R \colonequals  \{(q,r)\in \R^{c+d}: q\in Q, r\in R\}.
	\]
	The Cartesian product of two polytopes is a polytope and $\dim(Q\times R)=\dim(Q)+\dim(R)$. Furthermore every pair of nonempty faces $F_1\subset Q,F_2\subset R$ induces a nonempty face $F \colonequals F_1\times F_2$ of $Q\times R$ of dimension $\dim(F_1)+\dim(F_2)$. All nonempty faces of $Q\times R$ arise in this way. A \textbf{$d$-cube} is a polytope combinatorially isomorphic to the product of $d$ segments $\Delta_1$.

	\subsection{Gale diagrams}
	We now come to a central tool for our results.
	
	\begin{definition}
		Let $\A=\{a_1,\dots, a_n\}\subset \R^d$ be a point configuration affinely spanning $\R^d$, and let $M_{\textbf{1},\A}$ be the matrix where the $i$-th column is $(1,a_i)\in \R^{d+1}$. A \textbf{Gale transform} of $\A$ is an $n$-point configuration $\gale(\A)= \{b_1,\dots, b_n\}\subset \R^{n-d-1}$ such that the row span of $M_{\textbf{1},\A}$ is orthogonal to the row span of $M_{\gale(\A)}$.
	\end{definition}
	
	The importance of Gale transforms stems from the fact that $\F(\A)$ can be read directly from $\gale(\A)$. More precisely \cite[Chapter 3, Theorem 1]{bible} states that
	\begin{equation}\label{eq:gale}
	\{a_{i_!},\dots,a_{i_k}\}\subset \A\text{ is a coface } \Longleftrightarrow \O\in\relint\left(\conv\left(\{b_{i_!},\dots,b_{i_k}\}\right)\right).
	\end{equation} 
	
	\begin{definition}
		We call two point configurations $\A_1$ and $\A_2$ \textbf{Gale equivalent} if there exists a bijection $\psi$ between them such that $\O\in\relint(\conv\left(Z\right))\Longleftrightarrow \O\in\relint(\conv\left(\psi(Z)\right))$ for any subset $Z\subset\A_1$. Any point configuration that is Gale equivalent to a Gale transform of $\A$ is called a \textbf{Gale diagram} of $\A$.
	\end{definition}
	
	Our definition follows the notation in \cite{bible}; what we call a Gale transform is called Gale diagram in other sources, see \cite[Notes Chapter 6]{ziegler}.
	
	If $\A$ is an $n$-point configuration in $\R^d$, then a Gale diagram of $\A$ is an $n$-point configuration in $\R^{n-d-1}$. Thus Gale diagrams are particularly helpful when the number of vertices is small. Gale diagrams have also found uses in algebraic geometry; see, e.g., \cite{eisenbud}.
	\begin{example}\label{ex:gale}
		Let $\A=\{(0,0,0),(2,0,0),(0,2,0),(2,2,0),(1,1,1),(1,1,-1)\}$ be a point configuration in $\R^3$. The convex hull of $\A$ is an octahedron. Since the matrices $M_\A$ and $M_{\gale(\A)}$ below have orthogonal row spaces, the point configuration \\ $\gale(\A)=\{(1,0),(-1,-1),(-1,-1),(1,0),(0,1),(0,1)\}$ is a Gale transform of $\A$.
		\[
		M_\A=\begin{bmatrix}1&1&1&1&1&1\\ 0&2&0&2&1&1 \\ 0&0&2&2&1&1 \\ 0&0&0&0&1&-1 \end{bmatrix}, \qquad
		M_{\gale(\A)}=\begin{bmatrix}1&-1&-1&1&0&0\\ 0&-1&-1&0&1&1 \end{bmatrix}
		\]
		 The set $\{a_1,a_2,a_5\}$ is a coface of \(\A\), and $\conv(b_1,b_2,b_5)$ is a triangle that contains the origin in the interior, illustrating Equation \eqref{eq:gale}.
		
		\begin{figure}[ht]
			
			\begin{tikzpicture}
			\node at (-5,2.4) {Point configuration $\A$};
			\node at (0.2,2.4) {A Gale transform $\A$};
			\node at (5.3,2.4) {A Gale diagram of $\A$};
			\begin{scope}[xshift=-7cm, yshift=-0.4cm,scale=1.5,
			x = {(1cm,0cm)}, y = {(.5cm,.3cm)}, z  = {(0cm,1.2cm)}]
			\coordinate (1) at (0,0,0);
			\coordinate (2) at (2,0,0);
			\coordinate (3) at (0,2,0);
			\coordinate (4) at (2,2,0);
			\coordinate (5) at (1,1,1);
			\coordinate (6) at (1,1,-1);
			
			\draw[thick] (1)--(2)--(5)--(1)--(6)--(2)--(4)--(5) (4)--(6);
			\draw[dashed] (5)--(3) (4)--(3) (3)--(1) (3)--(6);
			
			\foreach \p in {1,...,6}
			\node[shape=circle,draw,thick,fill=white] at (\p) {\p};
			
			\end{scope}
			
			\begin{scope}[xshift=0cm,scale=1.5]
			\coordinate (1) at (1,0);
			\coordinate (2) at (-1,-1);
			\coordinate (3) at (-1,-1);
			\coordinate (4) at (1,0);
			\coordinate (5) at (0,1);
			\coordinate (6) at (0,1);
			
			\draw[dotted] (1)--(2)--(5)--(1)--(6)--(2)--(4)--(5) (4)--(6);
			
			\node[shape=circle,draw,thick,fill=white] at (1) {1,4};
			\node[shape=circle,draw,thick,fill=white] at (2) {2,3};
			\node[shape=circle,draw,thick,fill=white] at (5) {5,6};
			
			\draw[fill=black] (0,0) circle [radius=.03];
			\node[above right] at (0,0) {$O$};

			\end{scope}
			
			\begin{scope}[xshift=4.6cm,scale=1.5]
			\coordinate (1) at (1,0);
			\coordinate (2) at (-0.9,-1.2);
			\coordinate (3) at (-1.2,-0.9);
			\coordinate (4) at (1.2,0.4);
			\coordinate (5) at (0,1.3);
			\coordinate (6) at (-.4,1);
			
			\draw[dotted] (1)--(2)--(5)--(1);
			\draw[dotted] (1)--(3)--(6)--(1);
			\draw[dotted] (4)--(2)--(6)--(4);
			\draw[dotted] (4)--(3)--(5)--(4);
			
			\foreach \p in {1,...,6}
			\node[shape=circle,draw,thick,fill=white] at (\p) {\p};

			\draw[fill=black] (0,0) circle [radius=.03];
			\node[above right] at (0,0) {$O$};

			\end{scope}
			
			\end{tikzpicture}
			\caption{The main objects of interest in Example \ref{ex:gale}.}
			\label{fig:basic_example}
		\end{figure}
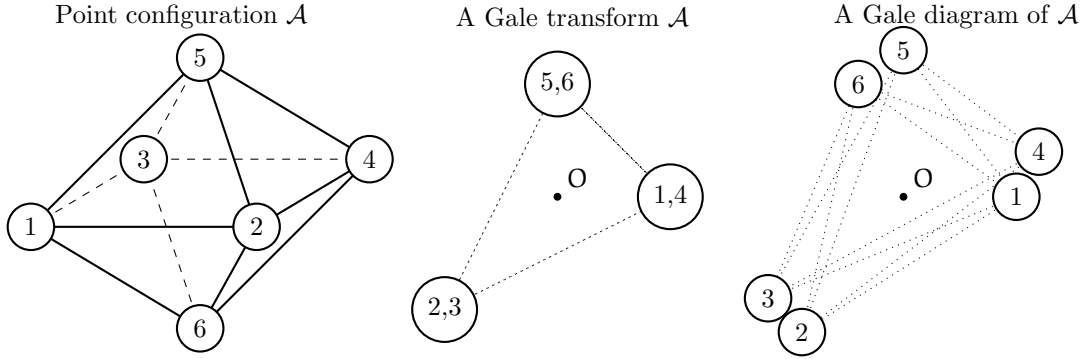
		
	\end{example}
	
	\begin{remark}
		We emphasize that our definition of Gale equivalent point configurations is weaker than isomorphism between the oriented matroids of the same point sets. In Example~\ref{ex:gale} the Gale diagram of \(\A\) shown at the right of Figure~\ref{fig:basic_example} has a different oriented matroid than the Gale transform of \(\A\) shown in the center of the same figure, since the pairs of points $(1,4),(2,3),(5,6)$ do not coincide. 
	\end{remark}	
	
	\section{Minkowski Summands}\label{sec:type}
	
	Let $Q,R \subset \R^n$ be two polytopes. We define their \textbf{Minkowski sum} to be
	\[
	Q+R \colonequals  \{q+r: q\in Q, r\in R\}.
	\]
	The Minkoswki sum of two polytopes is a polytope. We call $Q$ a weak Minkowski summand of $P$, denoted $Q\preceq P$, if there exist a polytope $R$ and a positive scalar $\lambda$ so that $Q+R=\lambda P$. 
	
	\begin{example}
		In Figure \ref{fig:sumex} we depict two polygons and their Minkowski sum. 
		\begin{figure}[ht]
			
			\begin{tikzpicture}
			
			\draw (1,1)--(2,3)--(3,2)--cycle;
			\node[right] at (2,2) {$Q$};
			
			\draw (1,5)--(3,3)--(4,4)--cycle;
			\node[right] at (3,4) {$R$};
			
			\foreach \x in {1,2,...,10}
			\foreach \y in {1,2,...,5}
			\draw[fill=black] (\x,\y) circle (.05);
			
			\begin{scope}[xshift=5cm, yshift=1cm]
			\draw (2,0)--(0,2)--(1,4)--(4,3)--(5,2)--(4,1)--cycle;
			\node[below] at (2,2) {$Q+R$};
			\end{scope}
			
			\end{tikzpicture}
			\caption{A Minkowski sum of two polytopes.}\label{fig:sumex}
		\end{figure}
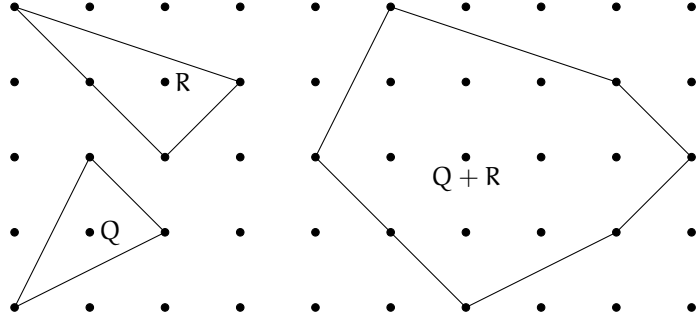
		
	\end{example}
	
	Given a polytope $P$, we are interested in the set of all its Minkowski summands. 
	
	\begin{definition} Let $P$ be a polytope. Let $\sim_T$ be the equivalence relation \(P_1 \sim_T P_2 \iff P_1 = P_2 + \vec{v}\) for some vector \(\vec{v}\). Let \(\sim_{D+T}\) be the equivalence relation \(P_1 \sim_{D+T} P_2 \iff P_1 = \lambda P_2 + \vec{v} \) for some scalar \(\lambda\) and some vector \(\vec{v}\). We define

		\begin{align}
		\Mink(P) &\colonequals \{Q\text{ a polytope}: Q\preceq P \},\\
		\TMink(P)&\colonequals \Mink(P)/\sim_T,\\
		\DMink(P)&\colonequals \Mink(P)/\sim_{D+T}.
		\end{align}
		
	\end{definition}
	
	In the following sections we will study two different ways to parameterize these sets as  polyhedral sets. For this we need some classical results characterizing (weak) Minkowski summands. We use the characterization given in \cite[Chapter 15]{grunb}, written by G. Shephard, which we restate in the form we need for the present paper.
	
	\begin{theorem}[Shephard \cite{grunb}]
		Let $P=\{x\in\R^d: Ux\leq z\}$ be an irredundant inequality description for a polytope with $m$ facets. For any polytope $Q\subset \R^d$ the following are equivalent.
		\begin{enumerate}\label{thm:wms}
			\item[(i)] $Q$ is a weak Minkoswki summand of $P$.
			\item[(ii)] There exists a map $\phi:\F_0(P)\to \F_0(Q)$ such that for $v_i,v_j\in \F_0(P)$ with $\{v_i,v_j\}\in \F_1(P)$ we have $\phi(v_i)-\phi(v_j) = \lambda_{i,j} (v_i-v_j)$, for some $\lambda_{i,j}\in\R_{\geq 0}$.
			\item[(iii)] There exists $\eta\in\R^m$ such that $Q=\{x\in\R^d: Ux\leq \eta\}$ and for any subset of rows \(S\) such that the linear system $\{\langle u_i,x\rangle=z_i, \forall i\in S\}$ defines a vertex of $P$, the linear system $\{\langle u_i,x\rangle=\eta_i, \forall i\in S\}$ defines a vertex in $Q$.
		\end{enumerate}
	\end{theorem}
	\begin{proof}
		This is essentially in \cite[Theorem 2, Chapter 15]{grunb} and its preliminaries. Condition $(iii)$ is taken as the definition of deformations in \cite{nested}. In \cite[Section 2]{mcmullen} (see also \cite[Appendix]{nested}), this condition is shown to be equivalent to the condition that the normal fan of $Q$ is a coarsening of the normal fan of $P$, which in turn is equivalent to \cite[Equation (3), Chapter 15]{grunb}. Condition $(ii)$ is \cite[Equation (2), Chapter 15]{grunb} without the bound $\lambda\leq1$ since we are dealing with weak Minkowski summands. Finally, the proof of \cite[Theorem 2, Chapter 15]{grunb} shows all these three conditions are equivalent. 
	\end{proof}
	
	\begin{remark}\label{rem:unique}
		In the description given in Theorem \ref{thm:wms}(iii) the vector $\eta$ is unique. In fact we have $\eta_i=\max_{x\in Q}\langle u_i,x\rangle$.
	\end{remark}
	
	\section{Parametrizing \(\TMink(P)\): Minkowski weights}\label{sec:polygons}

	We briefly describe the theory of 1-Minkowski weights. This framework allows us to study some classical polytopes, including the regular dodecahedron and polygons.
	
	\begin{definition}
		A \textbf{1-Minkowski weight} on $P$ is a function $\omega: \F_1(P)\to\R$ such that for each \(F \in \F_2(P)\) choosing a cyclic orientation \(v_E\)  of its edge vectors gives
		\begin{equation}\label{eq:balanced}
		\sum_{E\in F} v_E \cdot\omega(E) = \O.
		\end{equation}
		Equation \eqref{eq:balanced} is called the \textbf{balancing condition}.
		The set of all 1-Minkowski weights on $P$ is denoted $\Omega_1(P)$. See \cite{weights} for general information about Minkowski weights.
	\end{definition}
	
	\begin{example}\label{ex:dodeweigths}
		If $P$ is the regular dodecahedron, then via elementary geometry the balancing condition for each pentagonal face with edges $\{E_0,E_1,E_2,E_3,E_4\}$ is equal to 		$$\omega(E_0)+\omega(E_1)e^{\frac{2\pi}{5}i}+\omega(E_2)e^{\frac{4\pi}{5}i}+\omega(E_3)e^{\frac{6\pi}{5}i}+\omega(E_4)e^{\frac{8\pi}{5}i}=0,$$ 
		as complex numbers.
		\begin{figure}[ht]
			\begin{tikzpicture}

\begin{scope}[xshift=-3cm, scale=0.4]

\node at (0,-9) {A Minkowski weight vanishing on 10 edges.};

    \foreach \p in {1,2,...,5}
    	\draw[fill=black] (\p*72:8) circle (.05);	
    \foreach \p in {1,2,...,5}
    	\draw[fill=black] (\p*72:6) circle (.05);
    \foreach \p in {1,2,...,5}
    	\draw[fill=black] (\p*72+36:4) circle (.05);
    \foreach \p in {1,2,...,5}
    	\draw[fill=black] (\p*72+36:2) circle (.05);
	
	
     \foreach \p in {1,2,...,5}
    	\draw (\p*72:8)--(\p*72+72:8) (\p*72:8)--(\p*72:6) (\p*72:6)--(\p*72+36:4) (\p*72:6)--(\p*72-36:4) (\p*72+36:4)--(\p*72+36:2) (\p*72+36:2)--(\p*72-36:2);

    \foreach \p in {1,2,...,5}
        \node[red] at (\p*72+36:6.8) {1};
    \foreach \p in {1,2,...,5}
        \node[red] at (\p*72+4:7) {0};
    \foreach \p in {1,2,...,10}
        \node[red] at (\p*36+18:5) {$\phi$};
    \foreach \p in {1,2,...,5}
        \node[red] at (\p*72+43:3) {0};
    \foreach \p in {1,2,...,5}
        \node[red] at (\p*72:1.3) {1};      
   
\begin{scope}[xshift=-5cm,yshift=7cm, scale=2]
 \node at (0,0) {$\phi=\frac{\sqrt{5}-1}{2}$};
 \draw (-1,-.5)--(-1,.5)--(1,.5)--(1,-.5)--cycle;
\end{scope}
\end{scope}

\begin{scope}[xshift=3.5cm, scale=0.4]
\node at (0,-9) {A Minkowski weight vanishing};
\node at (0,-9.8) {on exactly one edge.};

    \foreach \p in {1,2,...,5}
    	\draw[fill=black] (\p*72:8) circle (.05);	
    \foreach \p in {1,2,...,5}
    	\draw[fill=black] (\p*72:6) circle (.05);
    \foreach \p in {1,2,...,5}
    	\draw[fill=black] (\p*72+36:4) circle (.05);
    \foreach \p in {1,2,...,5}
    	\draw[fill=black] (\p*72+36:2) circle (.05);
	
	
     \foreach \p in {1,2,...,5}
    	\draw (\p*72:8)--(\p*72+72:8) (\p*72:8)--(\p*72:6) (\p*72:6)--(\p*72+36:4) (\p*72:6)--(\p*72-36:4) (\p*72+36:4)--(\p*72+36:2) (\p*72+36:2)--(\p*72-36:2);

    \foreach \p in {1,2,...,5}
        \node[red] at (\p*72+36:6.8) {1};
    \foreach \p in {1,2,...,4}
        \node[red] at (\p*72+4:7) {1};
        
    \node[red] at (5*72+4:7) {$\gamma$};
    \foreach \p in {2,3,...,7}
        \node[red] at (\p*36+18:5) {1};
        
    \node[red] at (1*36+18:5) {$\gamma$}; \node[red] at (8*36+18:5) {$\gamma$}; \node[red] at (9*36+18:5) {$\beta$}; \node[red] at (10*36+18:5) {$\beta$};
    \foreach \p in {1,2,3}
        \node[red] at (\p*72+43:3) {$\gamma$};     
    \foreach \p in {4,5}
        \node[red] at (\p*72+43:3) {$\alpha$};
        
    \foreach \p in {1,4}
        \node[red] at (\p*72:1.3) {$\alpha$};
    \foreach \p in {2,3}
        \node[red] at (\p*72:1.3) {$\beta$};  
     \node[red] at (5*72:1.3) {0};  
   
\begin{scope}[xshift=5.7cm,yshift=6.5cm]
 \node[right] at (0,1.5) {$\alpha=\frac{3+\sqrt{5}}{2}$};
 \node[right] at (0,0) {$\beta=\frac{1}{2}$};
 \node[right] at (0,-1.5) {$\gamma=\frac{5+\sqrt{5}}{4}$}; 
 \draw (0,-2.5)--(0,2.5)--(4.2,2.5)--(4.2,-2.5)--cycle;
\end{scope}

\end{scope}
\end{tikzpicture}
			\caption{Two Minkowski weights on the regular dodecahedron. }
			\label{fig:dode}
		\end{figure}
	\end{example}
	
	\begin{definition}\label{def:type}
		Let $P$ be a polytope. We define
		\begin{align}
		\typecone(P) &\colonequals \left\{\omega\in\Omega_1(P): \omega(E)\geq 0, \forall E\in\F_1(P)\right\}, \label{eq:typecone}\\
		\type(P) &\colonequals \left\{\omega\in\typecone(P): \sum_{E\in\F_1(P)}\omega(E)=f_1(P)\right\}. \label{eq:type}
		\end{align}
		
		The \textbf{type cone} is the pointed polyhedral cone $\typecone(P)$.  We note that $\typecone(P)$ is a cone over the \textbf{type polytope} $\type(P)$, so they easily determine one another.
	\end{definition}
	
	The polyhedron $\typecone(P)$ parametrizes $\TMink(P)$.
	Indeed, Theorem \ref{thm:wms}(ii) guarantees the existence of a $1$-Minkowski weight for \(Q\) and conversely, \cite[Theorem 15.5]{faces} describes how to reconstruct $Q$ from a \(1\)-Minkowski weight, up to translation.
	The polytope $\type(P)$ parametrizes the set $\DMink(P)$.
	
	\begin{remark}
		We clarify a minor technical discrepency; in \cite{mcmullen}, ``type cone'' refers to the \emph{interior} of what we have defined as the \(\typecone(P)\), since every polytope in the interior of $\typecone(P)$, corresponding to stricly positive 1-Minkowski weigths, has the same combinatorial type. 
	\end{remark}
	
	\begin{definition}\label{def:vanishing}
		Let $P$ be a polytope and $S\subset\F_1(P)$ a subset of its edges. If there exists $\omega\in\type(P)$ such that $\omega(E)=0$ if and only if $E\in S$, then $S$ is a \textbf{vanishing} set. The faces of the type cone are in bijection with vanishing sets of edges. Translating this into \(\TMink(P)\), a vanishing set \(S\) corresponds to the set of Minkowski summands of \(P\) whose edges are in the complement of $S$.
	\end{definition}
	
	When $P$ is a simple $d$-dimensional polytope we have 
	\begin{equation}\label{eq:dim}
	dim \left(\typecone(P)\right) = f_{d-1}(P)-d,
	\end{equation}by \cite[Theorem 11]{mcmullen}. Hence $\dim \left(\type(P)\right) = f_{d-1}(P)-d-1.$ The dimension of $\typecone(P)$ is hard to compute in general. When $\dim(\typecone(P))=1$ we say that $P$ is an indecomposable polytope, since its only weak Minkowski summands are, up to translation, dilations of $P$. 
	
	\begin{remark}\label{rem:faces}
		Let $P$ be a polytope and $K=\typecone(P)$.
		The faces of $K$ are the type cones of the weak Minkowski summands of $P$ \cite[Theorem 7]{mcmullen}. Thus the rays of $K$ correspond to \emph{indecomposable} weak Minkowski summands. By Equation \eqref{eq:typecone}, $K$ has at most $f_1(P)$ facets. However, it can happen that some inequalities $\omega(E)\geq 0$ are not facet-defining, see Example~\ref{ex:polygons}. 
	\end{remark}
	
	\begin{example}[Regular Cubes]
		
		Let $P=[0,1]^d$ be a regular d-cube. Here the balancing condition is equivalent to all edges in each parallel class having same weight. Thus \(\typecone(P)\) can be parameterized by these $d$ weights, all of which must be nonnegative, so $\typecone(P)\cong \R^d_{\geq 0}$, the positive orthant.
	\end{example}
	
	\begin{example}[Regular Dodecahedron]
		
		Let $P$ be the regular dodecahedron. We can compute $\typecone(P)$ using \texttt{sage} \cite{sagemath}. 
		Since $P$ is a simple 3-polytope with 12 facets, we have $\dim(\typecone(P))=9$. Its $f$-vector is $$(1, 278, 2340, 6616, 8812, 6105, 2192, 375, 30, 1).$$ 
		The group of symmetries $G$ of the regular dodecahedron, the Coxeter group $H_3$, acts on $\type(P)$.
		The 30 facets correspond to the weights that vanish on a single edge, as shown on the right hand side of Figure \ref{fig:dode}. All of the facets are in the same $G$-orbit. 
		There are a total of $\binom{30}{2}=435$ pairs of edges.
		The $60$ adjacent pairs form a $G$-orbit and none of them are vanishing sets. This implies that all of the remaining $375$ pairs of edges are all vanishing sets. These pairs split into 7 orbits.
		
		The weight displayed on the left hand side of Figure \ref{fig:dode} describes a pentagonal antiprism, which is a vertex of $\type(P)$. There are 6 such antiprisms in the $G$-orbit corresponding to the 6 pairs of opposite faces, these account for only 6 of the 278 vertices of $\type(P)$.
	\end{example}

	\begin{example}[Polygons]\label{ex:polygons}
	    Let $P\subset \R^2$ be a $n$-gon.
	    Its edge vectors, oriented cyclically, give a vector configuration $\V(P)=\{v_1,\dots,v_n\}\subset \R^2$. 
	    In this case there is only one balancing condition so the type polytope $\type(P)$ is given by all vectors $(\omega_i)_{i\in[n]}\in\R^n$ satisfying
	    \begin{equation}\label{eq:poly}
	        \vec{v}_1\omega_1+\dots+\vec{v}_n\omega_n = 0,\quad \omega_1+\dots+\omega_n = n,\quad \text{and }\omega_i\geq 0\text{ for }i\in[n].
	    \end{equation}
	    This defines a $(n-3)$-dimensional polytope with at most $n$ facets.
	    
	    If $P$ is a pentagon, its type polytope can be a triangle, quadrilateral or a pentagon, see for example \cite[Figure 4]{bavard} for the triangle and quadrilateral case.
	    So the combinatorial type of $P$ does not determine the combinatorial type of $\type(P)$.
	    
	    Using Definition~\ref{def:polar}, we define \(\A(P)\) as the polar point configuration of system \eqref{eq:poly} consisting of $n$ points, one for each edge.
	    By the description of the faces of type polytopes in Definition \ref{def:vanishing}, it follows that $\V(P)$ is a Gale diagram for $\A(P)$.
	    Conversely, any set of $n$ nonzero distinct vectors in $\R^2$ having sum equal to zero are the edges of a polygon.
	    This shows that for every $(n-3)$-dimensional polytope $Q$ with at most $n$ facets there exists a polygon $P$ such that $Q$ and $\type(P)$ are combinatorially isomorphic.
	    
	    Notice that if we are given $n$ vectors in $\R^2$ having sum equal to zero together with an specified labeling using the set $[n]$, they may not be the edges of a convex polygon in the same order, but rather of a non-convex polygon.
	    See \cite{bavard} for more details about type cones of (possibly non-convex) polygons and the relation with hyperbolic geometry.
	\end{example}

	\section{Parametrizing $\Mink(P)$: Intersections in the Gale diagram}\label{sec:cubes}
	
	In \cite{mcmullen} McMullen gave a different technique to analyze type polytopes. In this section, we first discuss this technique and then apply it to compute the type polytope of the product of simplices.
	
	\begin{theorem}[McMullen \cite{mcmullen}]\label{thm:mcmullen}
		Let $P$ be a polytope, $\A=\{a_1,\dots,a_m\}$ be the vertex set of its polar $P^\circ$, and $\gale(A)=\{b_1,\dots,b_m\}$ be a Gale transform for $\A$. Then \[\type(P)\cong \bigcap_S \conv\{b_i: b_i\in S\},\] where the intersection is over all cofacets $S$ of $\A$.
	\end{theorem}
	
	Theorem~\ref{thm:mcmullen} follows from the results of \cite{mcmullen}, see in particular his comments on Page 88 at the end of Section 2. Since it is not explicitly stated in the source in the form we need, we sketch the main ideas of the proof.
	A detailed proof can also be found in \cite[Section 1]{fillastre}.

	\begin{proof}
		
		Let $P\subset \R^d$ be given by $P=\{x\in\R^d: Ux\leq z\}$ where $U$ is an $m\times d$ matrix that we identify with its set of rows $\{u_1,\dots, u_m\}$, and $z=(z_1,\dots,z_m)\in \R^m$. We assume each inequality is facet defining. For every vertex $v$ of $P$ we let $S_v:=\{u_i\in U: \langle u_i, v\rangle = z_i\}\subset U$, in other words the set of facets of $P$ that contain $v$.
		
		For each element $\eta\in \R^m$ we consider the (possibly empty) polytope $P_U(\eta) \colonequals \{x\in\R^d: Ux\leq \eta\}$. By Theorem~\ref{thm:wms}(iii), and Remark~\ref{rem:unique}, each element $Q\in\Mink(P)$ is represented as $Q=P_U(\eta)$ for a unique $\eta$, so that $\Mink(P)$ can be identified as a subset of $\R^m$.
		
		We want to consider the elements of $\Mink(P)$ up to translation. For any $w\in\R^d$ we have $P_U(\eta)+w=P_U(\eta+Uw)$. Let $\hU$ be an $(m-d)\times m$ matrix such that $\hU U= 0$ and consider the linear map $\phi:\R^d\longrightarrow \R^{m-d}$ given by matrix multiplication by $\hU$. Since $\phi(\eta+Uw)=\phi(\eta)$, $\phi$ maps the whole translation class of \(P_U(\eta)\) on to the same element. Thus we can identify $\TMink(P)$ with $\phi(\Mink(P))\subset \R^{m-d}$.
		
		By Theorem~\ref{thm:wms}(iii) we have that $Q\in\Mink(P)$ if and only if for every vertex $v$ of $P$ the following system has a (unique) solution: 
		\begin{align*}
		x&\in \R^d,\\
		\langle u_i,x\rangle &=\eta_i, \quad u_i\in S_v,\\
		\langle u_j,x\rangle &\leq \eta_i, \quad u_j\notin S_v.
		\end{align*}
		If there is a solution, up to translation, we can assume that the solution is $x=0$, so that $\eta_i=0$ whenever $u_i\in S_v$ and $\eta_i\geq 0$ whenever $u_i\notin S_v$. Let $\hu_1,\dots,\hu_m$ be the columns of $\hU$. So we have that $\eta$ determines a weak Minkowski summand if and only if $\phi(\eta):=\sum \hu_i\eta_i\in \cone(\hu_i: u_i\notin S_v)$ for every $v$. In other words
		\begin{equation}\label{eq:intersection}
		\phi(\Mink(P)) = \bigcap_{v} \cone(\hu_i: u_i\notin S_v),
		\end{equation}
		which we can identify with $\TMink(P)$ and in fact the cone described by Equation \eqref{eq:intersection} is linearly isomorphic to $\typecone(P)$.
		
		Finally to compute $\type(P)$ we restrict to an affine hyperplane $H$ intersecting\\ $\cone(\hu_1,\dots,\hu_m)$ and consider the points of intersection $b_i:= H \cap \R_+\hu_i$. The set $\{b_1,\dots,b_m\}$ is a Gale transform for $\{a_1,\dots,a_m\}$ where $a_i \colonequals z_i^{-1}u_i$ are the vertices of $P^\circ$. By polarity the complements of the sets $S_v$ considered in Equation~\eqref{eq:intersection} are exactly the cofacets of $P^\circ$.
	\end{proof}

	\subsection{Application: products of simplices}
	
	We begin by describing the face structure of arbitrary products of simplices. Let $(\Delta^0,\Delta^1,\dots,\Delta^k)$ be a list of $k+1$ simplices and let $d_i:=\dim(\Delta^i)\geq 1$ for each $i$. Throughout this section, we do not allow $0$-simplices to be factors in our product, since the cartesian product of $P$ with a point is isomorphic to $P$. We denote the vertex set of $\Delta^i$ by $\{v_{ij}\}_{0\leq j\leq d_i}$.
	Let $P=\prod_{i=0}^k \Delta^i$ be the Cartesian product of these \(k+1\) simplices. We have that $\dim(P)=D$, where $\sum_{i=0}^kd_i=D$. The following holds.
	\begin{enumerate}
		\item Vertices of $P$ are labeled with sequences $\{j_i\}_{0\leq i\leq k}$ where $0\leq j_i\leq d_i$ for each $i$. This labels the point \(v_{0j_0}\times\dots\times v_{kj_k}\).
		\item Facets of $P$ are labeled with pairs $(a,b)\in \{0,\dots,k\}\times\{0,\dots,d_a\}$. The facet \((a,b)\) contains the vertices labeled \(\{j_i\}_i\) such that \(j_a \neq b\). In other words, the facet labeled with $(a,b)$ is the set of all vertices of $P$ that do not have $v_{ab}$ as a factor.
		\item For a fixed vertex labeled $\{j_i\}_i$, the set of facets containing it are those labeled by $(a,b)$ such that $j_{a}\neq b$ for $0\leq a\leq k$.
	\end{enumerate}
	
	\begin{lemma}\label{lem:product}
		Let $P$ a product of $k+1$ simplices, then the cofacets of \(P^\circ\) are of size \(k+1\), and the vertices of $P^\circ$ can be colored with $k+1$ colors such that every cofacet contains a vertex of each color.
	\end{lemma}
	
	\begin{proof}
		The size of the cofacets follows directly from (3) above. Vertices are labeled by pairs $(a,b)$ as in (2) above and we assign vertex $p_{ab}$ to have color $a$. Under this coloring, each cofacet contains a vertex of each color.
	\end{proof}
	
	This structure motivates the following definition.

	\begin{definition}
		Let $\bd=(d_0,\dots,d_k)\in\N^{k+1}$ 
		and let $\rR(\bd)$ be a point configuration in $\R^k$  with a colored partition $\rR(\bd)=S_0\sqcup S_1\sqcup \dots \sqcup S_k$ such that each $S_i=\{p_{ij}\}_{j\in[d_i]}$ contains $d_i$ points of color $i$ for each $i=0,\dots,k$. A \textbf{rainbow subset} of \(\rR(\bd)\) is a subset $Z\subset\rR$ such that $|Z \cap S_i|=1$ for each $i=0,\dots,k$. A \textbf{rainbow simplex} is the convex hull of an affinely independent rainbow subset.
		The configuration $\rR(\bd)$ is a \textbf{rainbow configuration} if every rainbow subset is affinely independent and the intersection of all rainbow simplices is full dimensional. 
		
	\end{definition}
	
	\begin{remark}
	In the context of colored set partitions, there is a recent notion of generalized Gale transform studied in \cite{colorgale}.
	\end{remark}
	
	\begin{proposition}\label{prop:rainbow}
		Let $P\subset \R^d$ be combinatorially isomorphic to a product of $k+1$ simplices. Then every Gale transform $\G$ of $P^\circ$  is a rainbow configuration in $\R^{k}$.
	\end{proposition}	
	\begin{proof}
		A simple dimension count shows that $\G\subset\R^k$.
		
		In Lemma \ref{lem:product} we have already shown how to color the vertex set of $P^\circ$ and we can naturally carry this coloring to $\G$. By the same lemma, the cofacets are the rainbow subsets so following Equation \eqref{eq:gale} we have that every rainbow subset of \(\G\) contains the origin in its relative interior. 
		
		Suppose some rainbow subset \(S\) of \(\G\) were not affinely independent. Then some strict subset of \(S\) would contain the origin in its relative interior, and therefore correspond to a coface of \(P\). But the corresponding face of \(P\) would then strictly contain a facet, a contradiction. Therefore every rainbow subset \(S\) of \(\G\) is affinely independent.
		
		The rainbow subsets of \(\G\) each contain \(k+1\) affinely independent points, so the corresponding rainbow simplices are \(k\) dimensional and live in \(\mathbb{R}^k\). Therefore all rainbow simplices of \(\G\) are full dimensional. This shows that the intersection of all the rainbow simplices in \(\G\) is full dimensional, as each contain a neighborhood of the origin. We have therefore given \(\G\) a colored partition so that every rainbow subset is affinely independent, and so that the intersection of all rainbow simplices is full dimensional.
	\end{proof}
	
	\begin{example}
		Figure~\ref{fig:rain_config} shows an example of a rainbow configuration in $\R^2$. Note that the intersection of all rainbow triangles is itself a triangle. 
		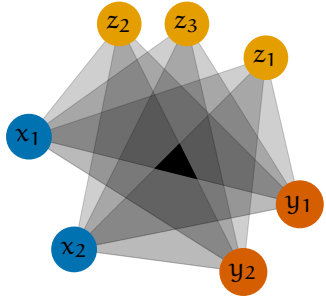
\begin{figure}[ht]
			\begin{tikzpicture}[scale=0.6]
			\usetikzlibrary{intersections}
			\node[coordinate] at (-3,1) (x1){};
			\node[coordinate] at (-2,-1.5) (x2){};
			\node[coordinate] at (3,-.5) (y1){};
			\node[coordinate] at (1.75,-2) (y2){};
			\node[coordinate] at (2.25,2.75) (z1){};
			\node[coordinate] at (-1,3.5) (z2){};
			\node[coordinate] at (0.5,3.5) (z3){};
			\draw[fill=black,opacity=0.1] (x1) -- (y1) -- (z1) -- cycle;
			\draw[fill=black,opacity=0.1] (x1) -- (y1) -- (z2) -- cycle;
			\draw[fill=black,opacity=0.1] (x1) -- (y1) -- (z3) -- cycle;
			\draw[fill=black,opacity=0.1] (x1) -- (y2) -- (z1) -- cycle;
			\draw[fill=black,opacity=0.1] (x1) -- (y2) -- (z2) -- cycle;
			\draw[fill=black,opacity=0.1] (x1) -- (y2) -- (z3) -- cycle;
			\draw[fill=black,opacity=0.1] (x2) -- (y1) -- (z1) -- cycle;
			\draw[fill=black,opacity=0.1] (x2) -- (y1) -- (z2) -- cycle;
			\draw[fill=black,opacity=0.1] (x2) -- (y1) -- (z3) -- cycle;
			\draw[fill=black,opacity=0.1] (x2) -- (y2) -- (z1) -- cycle;
			\draw[fill=black,opacity=0.1] (x2) -- (y2) -- (z2) -- cycle;
			\draw[fill=black,opacity=0.1] (x2) -- (y2) -- (z3) -- cycle;
			\node[circle,inner sep=2pt,fill=blue] at (x1){$x_1$};
			\node[circle,inner sep=2pt,fill=blue] at (x2){$x_2$};
			\node[circle,inner sep=2pt,fill=red] at (y1){$y_1$};
			\node[circle,inner sep=2pt,fill=red] at (y2){$y_2$};
			\node[circle,inner sep=2pt,fill=orange] at (z1){$z_1$};
			\node[circle,inner sep=2pt,fill=orange] at (z2){$z_2$};
			\node[circle,inner sep=2pt,fill=orange] at (z3){$z_3$};
			
			\path[name path=line 1] (z1) -- (x2);
			\path[name path=line 2] (z2) -- (y2);
			\path[name path=line 3] (y1) -- (x1);	
			\path [name intersections={of = line 1 and line 2}];
			\coordinate (c1)  at (intersection-1);
			\path [name intersections={of = line 1 and line 3}];
			\coordinate (c2)  at (intersection-1);
			\path [name intersections={of = line 2 and line 3}];
			\coordinate (c3)  at (intersection-1);
			\draw[fill=black,opacity=1] (c1) -- (c2) -- (c3) -- cycle;
			\end{tikzpicture}
			\caption{A rainbow configuration in the plane with three colors.}
			\label{fig:rain_config}
		\end{figure}
	\end{example}
	
	For the rest of the section let $\bd=(d_0,\dots,d_k)\in\N^{k+1}$ and let $\rR(\bd)$ be a rainbow configuration in $\R^k$ with point counts given by \(\bd\). 
	An \textbf{(affine) hyperplane} $H\subset\R^m$ is defined as $H=\{x\in \R^d: \langle y,x\rangle=b\}$, for some $y\in \R^d\backslash\{0\},b\in\R$.  Each affine hyperplane $H$ separates $\R^d\backslash H$ into two regions: $H_+=\{x\in \R^d: \langle y,x\rangle>b\}$ and $H_-=\{x\in \R^d: \langle y,x\rangle<b\}$. 
	We call an affine hyperplane $H$ \textbf{happy} with respect to the rainbow configuration \(\rR(\bd)\) if \(H\) is affinely spanned by points of $k$ distinct colors, every point of these \(k\) colors is in $H_{\leq 0} \colonequals  H\cup H_-$, and every point of the last color is in $H_+$. 
	
	\begin{proposition}\label{prop:atleast}
		Let $T$ be the intersection of all rainbow simplices of a rainbow configuration \(\rR(\bd)\) in \(\mathbb{R}^k\). Then the affine span of each facet of $T$ is a happy hyperplane.
	\end{proposition}
	
	\begin{proof}
		Since \(\rR(d)\) is a rainbow configuration, $T$ is a full dimensional polytope. Without loss of generality, fix some facet of \(T\). By the definition of $T$, the affine hull of this facet \(H\) contains a facet of some rainbow simplex. Since \(H\) contains the vertices of that facet, it must contain \(k\) points of distinct colors. Without loss of generality, assume that it contains the points \(\{p_1, \ldots, p_k\} \subset \rR(\bd)\) so that \(p_i\) is of color \(i\). We will show that $H$ is happy, having already showed the first condition.
		
		Let \(p_0\) and \(p_0'\) be points in \(\rR(\bd)\) of color \(0\). The rainbow simplices \(\conv\{p_0, p_1, \ldots, p_k\}\) and \(\conv\{p_0', p_1, \ldots, p_k\}\) must have full dimensional intersection. Since these simplices share a facet that lies in \(H\), \(p_0'\) and \(p_0\) must be on a common side of \(H\), say \(H^+\), showing the third condition of happiness.
		
		Without loss of generality, let \(p_k'\in S_k \subset \rR(\bd)\) be some point not of color \(0\).  Since \(H\) contains a facet of \(T\), \(H \cap \conv\{p_0, p_1, \ldots, p_{k-1}, p_k'\}\) must be \((k-1)\)-dimensional. As \(H\) contains the points \(p_1,\ldots,p_{k-1}\), the points \(p_0\) and \(p_k'\) must be on weakly different sides of \(H\). Therefore \(p_k'\) is in \(H_{\leq 0}\), showing the second condition of happiness.
		
		We therefore conclude that \(H\) is happy.
	\end{proof}
	
	\begin{proposition}\label{prop:atmost}
		Let \(\rR(\bd)\) be a rainbow configuration in \(\mathbb{R}^k\). Then there are at most $k+1$ happy hyperplanes, each one missing a distinct color. 
	\end{proposition}
	\begin{proof}
		Assume we have two happy hyperplanes $H_1,H_2$, each missing the color \(0\).
		Consider the polytope $Q_0 = \conv\{S_i: 1\leq i\leq k\},$ the convex hull of all points of \(\rR(\bd)\) not colored $0$. Each of the hyperplanes $H_1,H_2$ intersect \(Q_0\) on a facet, since they contain an affinely independent collection of \(k\) points in \(Q_0\), and all vertices of \(Q_0\) lie weakly on the same side of the \(H_i\). Let the corresponding facets of \(Q_0\) be \(F_1\) and \(F_2\). Pick any $p_0\in S_{0}$ and consider the polytope $Q'_0=\conv(Q_0,p_0)$. 
		From \cite[Lemma 4.3.2]{triangulations}, we see that $Q_0$ together with all pyramids to \(p_0\) over the facets of $Q_0$ whose supporting hyperplanes separate $p_0$ from \(Q_0\) (i.e., facets \emph{visible} from \(p_0\)) give a polyhedral subdivision of $Q'_0$. Since \(H_1\) and \(H_2\) are happy, \(F_1\) and \(F_2\) appear among these facets. As a consequence, the intersection of the pyramids $\textrm{pyr}_{p_0}(F_1)$ and $\textrm{pyr}_{p_0}(F_2)$ is not full dimensional. Each of these pyramids contains a rainbow simplex, contradicting the hypothesis that the intersection of all rainbow simplices is full dimensional. This finishes the proof.
	\end{proof}
	
	\begin{theorem}\label{thm:intersection}
		Let $T$ be the intersection of all rainbow simplices of some rainbow configuration \(\rR(\bd)\) in \(\mathbb{R}^k\). Then $T$ is a simplex.
	\end{theorem}
	
	\begin{proof}
		Each facet of $T$ spans a happy hyperplane by Proposition~\ref{prop:atleast}. By Propostion~\ref{prop:atmost}, $T$ has at most $k+1$ facets. Finally, since \(T\) is \(k\)-dimensional, \(T\) must be a simplex.
	\end{proof}

	We preface our final theorem with a comment on \(\type(P)\) for \(P\) a \(d\)-cube. In the parameterization of \(\type(P)\) by \(1\)-Minkowski weights using Equation \eqref{eq:typecone}, \(\type(P)\) is embedded in \(\mathbb{R}^{d2^{d-1}}\) given by exponentially many inequalities. It turns out that it has only $d+1$ facets.
	We now come to the crux of this paper. 
	
	\begin{theorem}\label{thm:main}
		For any $P$ combinatorially isomorphic to a product of $k+1$ simplices, $\type(P)$ is a simplex of dimension $k$. In particular, the type cone of any combinatorial cube is simplicial of the same dimension.
	\end{theorem}
	
	\begin{proof}
		
		Theorem~\ref{thm:mcmullen} and Proposition~\ref{prop:rainbow} together show that \(\type(P)\) is the intersection of all rainbow simplices in a rainbow configuration, and Theorem~\ref{thm:intersection} shows this intersection is a simplex.
	\end{proof}

	\bibliography{biblio}
	\bibliographystyle{amsplain}
	
	\vspace{.15in}
	
	{\sc Federico Castillo}, Pontificia Universidad Cat\'olica de Chile, federico.castillo@mat.uc.cl
	
	{\sc Joseph Doolittle}, TU Graz, jdoolittle@tugraz.at
	
	{\sc Bennet Goeckner}, University of Washington, goeckner@uw.edu
	
	{\sc Michael S. Ross}, Iowa State University, msross@iastate.edu
	
	{\sc Li Ying}, Vanderbilt University, li.ying@vanderbilt.edu

\end{document}